\documentclass[12pt]{article}

\usepackage{amssymb}

\newtheorem{thm}{Theorem}[section]
\newtheorem{defn}[thm]{Definition}
\newtheorem{prop}[thm]{Proposition}
\newtheorem{cor}[thm]{Corollary}
\newtheorem{lemma}[thm]{Lemma}
\newtheorem{rema}[thm]{Remark}

\newcommand{\halmos}{\rule{1ex}{1.4ex}}

\newcommand{\bea}{\begin{eqnarray}}
\newcommand{\eea}{\end{eqnarray}}
\newcommand{\nn}{\nonumber \\}
\newcommand{\be}{\begin {equation}}

\newcommand{\ee}{\end{equation}}

\newcommand{\cc}{\varphi}

 \newcommand{\res}{\mbox{\rm Res}}
\renewcommand{\hom}{\mbox{\rm Hom}}
 \newcommand{\pf}{{\it Proof.}\hspace{2ex}}
 \newcommand{\epf}{\hspace*{\fill}\mbox{$\halmos$}}
 \newcommand{\epfv}{\hspace*{\fill}\mbox{$\halmos$}\vspace{1em}}
 \newcommand{\epfe}{\hspace{2em}\halmos}

\title{ {\bf Intertwining operator superalgebras 
and vertex tensor categories
for superconformal algebras, II} }
\date{}
\author{Yi-Zhi Huang and Antun Milas}
\begin{document}

\bibliographystyle{alpha}
\maketitle

\begin{abstract}
We  construct the intertwining operator superalgebras 
and vertex tensor
categories for the $N=2$ superconformal unitary minimal models and 
other related models. 
\end{abstract}

\renewcommand{\theequation}{\thesection.\arabic{equation}}
\renewcommand{\thethm}{\thesection.\arabic{thm}}
\setcounter{equation}{0}
\setcounter{thm}{0}
\setcounter{section}{-1}

\section{Introduction}

It has been known that the $N=2$ Neveu-Schwarz superalgebra is one of
the most important algebraic objects realized in superstring theory.
The $N=2$ superconformal field theories constructed {}from its discrete
unitary representations of central charge $c<3$ are among the
so-called ``minimal models.'' In the physics literature, there have
been many conjectural connections among Calabi-Yau manifolds,
Landau-Ginzburg models and these $N=2$ unitary minimal models. In fact,
the physical construction of mirror manifolds \cite{GP} 
used the conjectured
relations  \cite{G1}
\cite{G2} between certain particular Calabi-Yau manifolds 
and certain $N=2$ superconformal field theories (Gepner models)
constructed {}from unitary minimal models (see \cite{Gr} for a survey).  
To establish these
conjectures as mathematical theorems, it is necessary to construct the
$N=2$ unitary minimal models mathematically and to study their
structures in detail.

In the present paper, we apply the theory of intertwining operator algebras
developed by the first author in \cite{H2.5}, \cite{H4} and \cite{H6}
and the tensor product theory for modules
for a vertex operator algebra developed by Lepowsky and the first author
in \cite{HL1}--\cite{HL5}, \cite{HL6} and \cite{H1}
to construct the intertwining operator algebras and vertex tensor categories
for $N=2$ superconformal unitary minimal models. The main work in this
paper is to verify that the conditions to use the general theories  are 
satisfied for these models. The main technique used 
is the representation theory
of the $N=2$ Neveu-Schwarz algebras, which has been studied by many
physicists and mathematicians, especially by 
Eholzer and Gaberdiel \cite{EG}, Feigin, Semikhatov, Sirota and Tipunin 
\cite{FSST} \cite{FST}, and  Adamovi\'c \cite{A} \cite{A2}.

The present paper is organized as follows: In Section 1, we recall 
the notion of $N=2$ superconformal vertex operator superalgebra.
In Section 2, we recall and prove
some basic results on representations
of unitary minimal $N=2$ superconformal vertex operator superalgebras
and on representations of $N=2$ superconformal 
vertex operator superalgebras in a much more general class.
 Section 3 is devoted to 
the proof of the convergence and extension properties for 
products of intertwining operators for unitary
minimal $N=2$ superconformal vertex operator superalgebras and 
for vertex operator superalgebras in the general class. 
Our main results on
the intertwining operator superalgebra structure and 
vertex tensor category structure are given in Section 4.

{\bf Acknowledgment}: The research of 
Y.-Z. H. is supported in part by
NSF grant DMS-9622961.

\renewcommand{\theequation}{\thesection.\arabic{equation}}
\renewcommand{\thethm}{\thesection.\arabic{thm}}
\setcounter{equation}{0}
\setcounter{thm}{0}

\section {$N=2$ superconformal vertex operator superalgebras}

In this section we recall the notion of $N=2$ superconformal vertex
operator algebra and basic properties of such an algebra. These
algebras have been studied extensively by physicists. The following 
precise version of the definition is {}from \cite{A} and \cite{HZ}:

\begin{defn}
{\rm An {\it $N=2$ superconformal vertex operator superalgebra} is a
vertex operator superalgebra $(V, Y, \mathbf{1}, \omega)$ together with
odd elements $\tau^{+}$, $\tau^{-}$ and an even element $\mu$ satisfying the
following axioms: Let
\begin{eqnarray*}
Y(\tau^{+}, x)
&=&\sum_{r\in \mathbb{Z}+\frac{1}{2}}G^{+}(r)x^{-r-3/2},\\
Y(\tau^{-}, x)
&=&\sum_{r\in \mathbb{Z}+\frac{1}{2}}G^{-}(r)x^{-r-3/2},\\
Y(\mu, x)&=&\sum_{n\in \mathbb{Z}}J(n)x^{-n-1}.
\end{eqnarray*}
Then $V$ is a direct sum of eigenspaces of $J(0)$ with integral
eigenvalues which modulo $2\mathbb{Z}$ give the $\mathbb{Z}_{2}$ grading
for the vertex operator superalgebra structure, and  the following 
$N=2$ Neveu-Schwarz relations hold: For $m, n\in \mathbb{Z}$, $r,s \in
\mathbb{Z}+\frac{1}{2}$
\begin{eqnarray*}
{[L(m), L(n)]}&=&(m-n)L(m+n)+\frac{c}{12}(m^{3}-m)\delta_{m+n, 0},\\
{[J(m), J(n)]}&=&\frac{c}{3}m\delta_{m+n, 0},\\
{[L(m), J(n)]}&=&-nJ_{m+n},\\
{[L(m), G^{\pm}(r)]}&=&\left(\frac{m}{2}-r\right)G^{\pm}(m+r),\\
{[J(m), G^{\pm}(r)]}&=&\pm
G^{\pm}(m+r),\\
{[G^{+}(r), G^{-}(s)]}&=& 
2L(r+s)+(r-s)J(r+s)+\frac{c}{3}(r^{2}-\frac{1}{4})\delta_{r+s, 0},\\
{[G^{\pm}(r), G^{\pm}(s)]}&=&0
\end{eqnarray*}
where $L(m)$, $m\in \mathbb{Z}$, are the Virasoro operators on $V$ and
$c$ is the central charge of $V$.

{\it Modules} and {\it intertwining operators} for an $N=1$ superconformal 
vertex operator superalgebra are modules and intertwining operators
for the underlying vertex operator superalgebra.}
\end{defn}

The $N=2$ superconformal vertex operator superalgebra defined above is
denoted by $(V, Y, \mathbf{1}, \tau^{+}, \tau^{-}, \mu)$ 
(without $\omega$ since 
$$\omega=L(-2)\mathbf{1}=\frac{1}{2}
[G^{+}(-\frac{3}{2}), G^{-}(-\frac{1}{2})]\mathbf{1}
+\frac{1}{2}J(-2)\mathbf{1})$$
or simply by $V$. 
Note that
a module $W$ for a vertex operator superalgebra 
(in particular the algebra itself) 
has a
$\mathbb{Z}_{2}$-grading called {\it sign} in addition to the
$\mathbb{C}$-grading by weights. We shall always use 
$W^{0}$ and $W^{1}$ to denote the even and odd subspaces of 
$W$. If $W$ is irreducible, there exists
$h\in \mathbb{C}$ such that $W=W^{0}\oplus W^{1}$ where
$W^{0}=\oplus_{n\in h+\mathbb{Z}}W_{(n)}$ and $W^{1}=\oplus_{n\in
h+\mathbb{Z}+1/2}W_{(n)}$ are the even and odd parts of $W$,
respectively. We shall always use the notation $|\cdot|$ to denote the
map {}from the union of the even and odd subspaces of a
vertex operator superalgebra or of 
a module for such an algebra to $\mathbb{Z}_{2}$ 
by taking the signs of elements in the union.

The notion of $N=2$ superconformal vertex operator superalgebra 
can be reformulated using odd formal variables. (In the 
$N=1$ case, this reformulation was given by Barron
\cite{B1} \cite{B2}.) 

As in Part I (\cite{HM}), for $l$ symbols
$\varphi_{1}, \dots, \varphi_{l}$, consider
the exterior algebra of
the vector space over $\mathbb{C}$ spanned by these symbols  
and denote this exterior algebra by $\mathbb{C}[\cc_{1}, \dots,
\cc_{l}]$. For any 
vector space $E$, we also have 
the vector space  
\begin{eqnarray*}
&E[\cc_{1}, \dots,
\cc_{l}],&\\
&E[x_{1}, \dots, x_{k}][\cc_{1}, \dots, \cc_{l}],&\\
&E[x_{1}, x_{1}^{-1}, 
\dots, x_{k}, x_{k}^{-1}][\cc_{1}, \dots, \cc_{l}],&\\
&E[[x_{1}, \dots, x_{k}]][\cc_{1}, \dots, \cc_{l}],&\\
&E[[x_{1}, x_{1}^{-1}, 
\dots, x_{k}, x_{k}^{-1}]][\cc_{1}, \dots, \cc_{l}],&\\
&E\{x_{1}, 
\dots, x_{k}\}[\cc_{1}, \dots, \cc_{l}]&
\end{eqnarray*}
and
$$E((x_{1}, 
\dots, x_{k}))[\cc_{1}, \dots, \cc_{l}].$$
If $E$ is a $\mathbb{Z}_{2}$-graded vector space, then there 
are  natural structures of modules over the ring 
$\mathbb{C}[x_{1}, \dots, x_{k}][\cc_{1}, \dots, \cc_{l}]$
on these spaces.

Let $(V, Y, \mathbf{1}, \tau^{+}, \tau^{-}, h)$ be an 
$N=2$ superconformal vertex operator superalgebra. Let
\begin{eqnarray*}
\tau_{1}&=&\frac{1}{\sqrt{2}}(\tau^{+}+\tau^{-}),\\
\tau_{2}&=&\frac{1}{\sqrt{-2}}(\tau^{+}-\tau^{-})
\end{eqnarray*}
and 
\begin{eqnarray*}
Y(\tau_{1}, x)&=&\sum_{r\in \mathbb{Z}+\frac{1}{2}}G_{1}(r)x^{-r-3/2},\\
Y(\tau_{2}, x)&=&\sum_{r\in \mathbb{Z}+\frac{1}{2}}G_{2}(r)x^{-r-3/2}.
\end{eqnarray*}
We define the {\it vertex operator map with odd variables}
\begin{eqnarray*}
Y: V\otimes V&\to &V((x))[\varphi_{1}, \varphi_{2}]\\
u\otimes v&\mapsto& Y(u, (x, \varphi_{1}, \varphi_{2})) v
\end{eqnarray*}
by
\begin{eqnarray*}
Y(u, (x, \varphi_{1}, \varphi_{2}))v&=&Y(u, x)v+\varphi_{1} 
Y(G_{1}(-1/2)u, x)v\nn
&&+ \varphi_{2}
Y(G_{2}(-1/2)u, x)v\nn
&&+\varphi_{1} \varphi_{2}Y(G_{1}(-1/2)
G_{2}(-1/2)u, x)v
\end{eqnarray*}
for $u, v\in V$. (We use the same notation $Y$ to denote 
the vertex operator map and the vertex operator map with odd
variables.)
In particular, we have:
\begin{eqnarray*}
Y(h,(x,\cc_1,\cc_2))&=&Y(h,x) -\sqrt{-1}\cc_1Y(\tau_2,x) \nn
&& -\sqrt{-1}\cc_2 Y(\tau_1,x)-2\sqrt{-1}\cc_1\cc_2 Y(\omega,x).
\end{eqnarray*}
Also, if we introduce 
\begin{eqnarray*}
\varphi^+&=&\frac{-\cc_1+\sqrt{-1}\cc_2}{\sqrt{2}}\\
\varphi^-&=&\frac{\cc_1+\sqrt{-1}\cc_2}{\sqrt{2}},
\end{eqnarray*}
then we can write
\begin{eqnarray*}
Y(h,(x,\cc_1,\cc_2))&=&Y(h,x)+\varphi^+Y(\tau^+,x) \nn
&&+\varphi^-Y(\tau^-,x)+2\varphi^+\varphi^-Y(\omega,x).
\end{eqnarray*}

We have:

\begin{prop}\label{odd-svoa}
The vertex operator map with odd variables satisfies the following 
properties:

\begin{enumerate}

\item The {\it vacuum property}: 
$$Y(\mathbf{1}, (x, \varphi_{1}, \varphi_{2}))=1$$
where $1$ on the right-hand side is the identity map on $V$.

\item The {\it creation property}: For any $v\in V$,
$$ Y(v, (x, \varphi_{1}, \varphi_{2}))\mathbf{1} \in V[[x]][\cc_1,\cc_2],$$
$$ \lim_{(x, \varphi_{1}, \varphi_{2}) 
\mapsto (0, 0, 0)} Y(v, (x, \varphi_{1}, \varphi_{2}))\mathbf{1}=v.$$

\item The {\it Jacobi identity}: In $(\mbox{\rm End}\ V)[[x_{0}, x_{0}^{-1}, 
x_{1}, x_{1}^{-1}, 
x_{2}, x_{2}^{-1}]][\cc_{1}, \cc_{2}, \psi_{1}, \psi_{2}]$, we have 
\bea
\lefteqn{x_0^{-1} \delta \left ( \frac {x_1-x_2-\cc_1
\psi_1-\cc_2
\psi_2}{x_0} \right ) Y(u,(x_1,\cc_1,\cc_2)) Y(v,(x_2,\psi_1,\psi_2))} \nn
&&\quad -(-1)^{|u||v|}x_0^{-1} \delta \left ( \frac {x_2-x_1+\cc_1
\psi_1+\cc_2
\psi_2}{-x_0} \right ) \cdot\nn
&&\quad\quad\quad\quad \cdot Y(v,(x_2,\psi_1,\psi_2)) Y(u,(x_1,\cc_1,\cc_2))
\nn
&& =x_2^{-1} \delta \left ( \frac {x_1-x_0-\cc_1
\psi_1-\cc_2
\psi_2}{x_2} \right
)\cdot\nn
&&\quad\quad\quad\quad \cdot
Y(Y(u,(x_0,\cc_1-\psi_1,\cc_2-\psi_2))v,(x_2,\psi_1,\psi_2)),\quad\quad\quad
\eea
for $u, v\in V$ which are either even or odd.

\item The {\it $G_{i}(-1/2)$-derivative property}: For any $v\in V$,
$i=1, 2$,
$$Y(G_{i}(-1/2) v, (x, \cc_{1}, \cc_{2}))=
\left(\frac{\partial}{\partial \cc_{i}}+
\cc_{i}\frac{\partial}{\partial x}\right)Y( v, (x,\cc_{1}, \cc_{2})).$$

\item The {\it $L(-1)$-derivative property}: For any $v\in V$, 
$$Y(L(-1) v, (x,\cc_{1}, \cc_{2}))=\frac{\partial}{\partial x}
Y(v, (x,\cc_{1}, \cc_{2})).$$

\item The {\it skew-symmetry}: For any $u, v\in V$ which are either 
even or odd,
\bea
&&\!\!\!\!\!\!\!\!\!\!Y(u,(x,\cc_1,\cc_2))v \nn
&&\!\!\!\!\!\!=(-1)^{|u||v|}e^{xL(-1)+\cc_1 G_1(-1/2)+\cc_2 G_2(-1/2) }
Y(v,(-x,
-\cc_1,-\cc_2))u.\hspace{1em}\halmos\nn
\eea
\end{enumerate}

\end{prop}

The proof of this result is straightforward and we omit it.

As in the $N=1$ case \cite{HM}, we can also reformulate 
the data and axioms for modules 
and intertwining operators for an $N=2$ 
superconformal vertex operator superalgebra using odd variables. 
Here we give the details of 
the corresponding reformulation of the data and axioms 
for intertwining operators. 

Let $W_{1}$, $W_{2}$ and $W_{3}$ be modules for an $N=2$ 
superconformal vertex operator superalgebra $V$ and $\mathcal{Y}$ an 
intertwining operator of type ${W_{3}\choose W_{1}W_{2}}$. We define the 
corresponding {\it intertwining operator map with odd variable}
\begin{eqnarray*}
\mathcal{Y}: W_{1}\otimes W_{2}&\to &W_{3}\{x\}[\varphi_1,\varphi_2]\\
w_{(1)}\otimes w_{(2)}&\mapsto& \mathcal{Y}(w_{(1)}, (x, \varphi_1,
\varphi_2)) w_{(2)}
\end{eqnarray*}
by
\begin{eqnarray*}
\mathcal{Y}(w_{(1)}, (x, \varphi_{1}, \varphi_{2})) w_{(2)}
&=&\mathcal{Y}(w_{(1)}, x) w_{(2)}+
\varphi_{1}\mathcal{Y}(G_{1}(-1/2)w_{(1)}, x) w_{(2)}\nn
&&+\varphi_{2}\mathcal{Y}(G_{2}(-1/2)w_{(1)}, x) w_{(2)}\nn
&&
+\varphi_{1}\varphi_{2}\mathcal{Y}(G_{1}(-1/2)G_{2}(-1/2)w_{(1)}, x)
 w_{(2)}
\end{eqnarray*}
for $u, v\in V$.
Then we have:

\begin{prop}
The intertwining operator map with odd variable satisfies the following 
properties:

\begin{enumerate}

\item The {\it Jacobi identity}: In $\hom(W_{1}\otimes W_{2}, W_{3})
\{x_{0},
x_{1},
x_{2}\}[\cc_{1}, \cc_{2},\psi_1, \psi_2]$, we have 
\bea
\lefteqn{x_0^{-1} \delta \left ( \frac {x_1-x_2-\cc_1
\psi_1-\cc_2
\psi_2}{x_0} \right ) Y(u,(x_1,\cc_1,\cc_2)) 
\mathcal{Y}(w_{(1)},(x_2,\psi_1,\psi_2))} \nn
&&\quad -(-1)^{|u||w_{(1)}|}x_0^{-1} 
\delta \left ( \frac {x_2-x_1+\cc_1
\psi_1+\cc_2
\psi_2}{-x_0} \right ) \cdot\nn
&&\quad\quad\quad\quad \cdot\mathcal{Y}(w_{(1)},(x_2,\psi_1,\psi_2)) 
Y(u,(x_1,\cc_1,\cc_2))
\nn
&& =x_2^{-1} \delta \left ( \frac {x_1-x_0-\cc_1
\psi_1-\cc_2
\psi_2}{x_2} \right
)\cdot\nn
&&\quad\quad\quad\quad \cdot
\mathcal{Y}(Y(u,(x_0,\cc_1-\psi_1,\cc_2-\psi_2))w_{(1)},
(x_2,\psi_1,\psi_2)),\quad\quad\quad\nn
\eea
for $u\in V$ and $w_{(1)}\in W_{1}$ which are either even or odd.

\item The {\it $G_{i}(-1/2)$-derivative property}: For any $v\in V$,
$i=1, 2$,
$$\mathcal{Y}(G_{i}(-1/2) v, (x, \cc_{1}, \cc_{2}))=
\left(\frac{\partial}{\partial \cc_{i}}+
\cc_{i}\frac{\partial}{\partial x}\right)\mathcal{Y}( v, (x,\cc_{1},
\cc_{2})).$$

\item The {\it $L(-1)$-derivative property}: For any $v\in V$, 
$$\mathcal{Y}(L(-1) w_{(1)}, (x,\cc_{1}, \cc_{2}))
=\frac{\partial}{\partial x}
\mathcal{Y}(w_{(1)}, (x, \cc_{1}, \cc_{2})).$$

\item The {\it skew-symmetry}: There is a linear
isomorphism
$$\Omega: \mathcal{V}^{W_{3}}_{W_{1} W_{2}} \to
\mathcal{V}^{W_{3}}_{W_{2} W_{1}}
$$
such that 
\begin{eqnarray*}
\lefteqn{\Omega (\mathcal{ Y})(w_{(1)}, (x,\cc_1,\cc_2))w_{(2)}} \nn
&&=(-1)^{|w_{(1)}||w_{(2)}|}e^{xL(-1)+\cc_{1}
G_{1}(-1/2) +\cc_2 G_{2}(-1/2)}\cdot\nn
&&\hspace{4em}\cdot\mathcal{ Y}(w_{(2)}, 
(e^{-\pi i}x,-\cc_1,-\cc_2)))w_{(1)}
\end{eqnarray*}
for $w_{(1)}\in W_{1}$ and $w_{(2)}\in W_{2}$ which are either even 
or odd. \epf

\end{enumerate}

\end{prop}

The proof of this result is similar to the proof of Proposition 
\ref{odd-svoa} and is omitted.

\renewcommand{\theequation}{\thesection.\arabic{equation}}
\renewcommand{\thethm}{\thesection.\arabic{thm}}
\setcounter{equation}{0}
\setcounter{thm}{0}

\section{Unitary minimal $N=2$ superconformal vertex operator superalgebras}

In this section, we recall the constructions and results on unitary minimal
$N=2$ superconformal vertex operator superalgebras and their
representations. New results needed in later
sections are also proved.  
We then introduce in this section a class of $N=2$ superconformal vertex
operator superalgebras and generalize most of the results for unitary minimal
$N=2$ superconformal vertex operator superalgebras to algebras in this
class. 

The {\it $N=2$ Neveu-Schwarz Lie superalgebra} 
is the vector space
$$\oplus_{n\in \mathbb{Z}}\mathbb{C}L_{n}\oplus
\oplus_{r\in \mathbb{Z}+\frac{1}{2}}\mathbb{C}G^{+}_{r}\oplus 
\oplus_{r\in \mathbb{Z}+\frac{1}{2}}\mathbb{C}G^{-}_{r}\oplus 
\oplus_{n\in \mathbb{Z}}\mathbb{C}J_{n}\oplus \mathbb{C}C$$
equipped with  the following 
$N=2$ Neveu-Schwarz relations:
\begin{eqnarray*}
{[L_{m}, L_{n}]}&=&(m-n)L_{m+n}+\frac{C}{12}
(m^{3}-m)\delta_{m+n, 0},\\
{[J_{m}, J_{n}]}&=&\frac{c}{3}m\delta_{m+n, 0},\\
{[L_{m}, J_{n}]}&=&-nJ_{m+n},\\
{[L_{m}, G^{\pm}_{r}]}&=&\left(\frac{m}{2}-r\right)G^{\pm}_{m+r},\\
{[J_{m}, G^{\pm}_{r}]}&=&\pm
G^{\pm}_{m+r},\\
{[G^{+}_{r}, G^{-}_{s}]}&=& 
2L_{r+s}+(r-s)J_{r+s}+\frac{C}{3}(r^{2}-\frac{1}{4})\delta_{r+s, 0},\\
{[G^{\pm}_{r}, G^{\pm}_{s}]}&=&0
\end{eqnarray*}
for $m, n\in \mathbb{Z}$, $r, s\in \mathbb{Z}+\frac{1}{2}$. 
For simplicity, we shall simply denote the $N=2$
Neveu-Schwarz Lie superalgebra by $\mathfrak{n}\mathfrak{s}(2)$ in this paper.

We now construct representations of the $N=2$ Neveu-Schwarz 
Lie superalgebra.
Consider the two subalgebras 
\begin{eqnarray*} 
\mathfrak{n}\mathfrak{s}^{+}(2)&=&\oplus_{n>0}\mathbb{C}L_{n}\oplus
\oplus_{r> 0}\mathbb{C}G^{+}_{r}\oplus 
\oplus_{r> 0}\mathbb{C}G^{-}_{r}\oplus 
\oplus_{n>0}\mathbb{C}J_{n},\\
\mathfrak{n}\mathfrak{s}^{-}(2)&=&\oplus_{n<0}\mathbb{C}L_{n}\oplus
\oplus_{r<0}\mathbb{C}G^{+}_{r}\oplus 
\oplus_{r<0}\mathbb{C}G^{-}_{r}\oplus 
\oplus_{n<0}\mathbb{C}J_{n}
\end{eqnarray*}
of $\mathfrak{n}\mathfrak{s}(2)$. 
Let $U(\cdot)$ be the functor {}from the 
category of Lie superalgebras to
the category of associative algebras obtained by taking the universal 
enveloping algebras of  Lie superalgebras. 
For any representation of $\mathfrak{n}\mathfrak{s}(2)$, 
we shall use 
$L(n)$, $n\in \mathbb{Z}$,  $G^{\pm}(r)$, $r\in \mathbb{Z}+\frac{1}{2}$,
and $J(n)$, $n\in \mathbb{Z}$,
to denote the representation images of $L_{n}$,
$G^{\pm}(r)$ and $J_{n}$.

For any $c, h, q\in \mathbb{C}$, 
the Verma module $M_{\mathfrak{ns}(2)}(c, h, q)$ 
for $\mathfrak{n}\mathfrak{s}(2)$ 
is a free $U(\mathfrak{n}\mathfrak{s}^{-}(2))$-module generated by
$\mathbf{1}_{c, h, q}$ such that
\begin{eqnarray*}
\mathfrak{n}\mathfrak{s}^{+}(2)\mathbf{1}_{c, h, q}&=&0,\\
L(0)\mathbf{1}_{c, h, q}&=&h\mathbf{1}_{c, h, q},\\
C\mathbf{1}_{c, h, q}&=&c\mathbf{1}_{c, h, q},\\
J(0)\mathbf{1}_{c, h, q}&=&q\mathbf{1}_{c, h, q}.
\end{eqnarray*}
There exists a unique maximal proper submodule $J_{\mathfrak{ns}(2)}(c,
h, q)$
of $M_{\mathfrak{ns}(2)}(c, h, q)$. It 
is easy to see that when $c\ne 0$, 
$\mathbf{1}_{c, 0, 0}$, $G^{\pm}(-3/2)\mathbf{1}_{c, 0, 0}$ and
$L(-2)\mathbf{1}_{c, 0, 0}$ are not in 
$J_{\mathfrak{n}\mathfrak{s}(2)}(c, 0, 0)$.
Let 
$$
L_{\mathfrak{ns}(2)}(c, h, q)=M_{\mathfrak{ns}(2)}(c, h, q)
/J_{\mathfrak{ns}(2)}(c, h, q)
$$
and 
$$V_{\mathfrak{ns}(2)}(c, 0, 0)=M_{\mathfrak{ns}(2)}(c, 0, 0)/\langle 
G^{+}(-1/2)\mathbf{1}_{c, 0, 0}, G^{-}(-1/2)\mathbf{1}_{c, 0, 0}\rangle$$
where $\langle 
G^{+}(-1/2)\mathbf{1}_{c, 0, 0}, G^{-}(-1/2)\mathbf{1}_{c, 0, 0}\rangle$ 
is the submodule of 
$M_{\mathfrak{ns}(2)}(c, 0, 0)$ generated by $G^{\pm}(-1/2)
\mathbf{1}_{c, 0, 0}$.
Then $L_{\mathfrak{ns}(2)}(c, 0, 0)$ and $V_{\mathfrak{ns}(2)}(c, 0, 0)$
have the structures of vertex operator superalgebras with
the vacuum $\mathbf{1}_{c, 0, 0}$, the Neveu-Schwarz elements
$G^{\pm}(-3/2)\mathbf{1}_{c, 0, 0}$ and the Virasoro element
$L(-2)\mathbf{1}_{c, 0, 0}$  (see
\cite{A}). 

In \cite{EG}, Eholzer and Gaberdiel showed, among other things, that
among vertex operator superalgebras of the form
$L_{\mathfrak{ns}(2)}(c, 0, 0)$, the only ones having finitely many
irreducible modules are the ``unitary'' ones $L_{\mathfrak{ns}(2)}(c_{m}, 0,
0)$ for nonnegative integers $m$, where $c_{m}=\frac{3m}{m+2}$.  The
following result was proved by Adamovi\'{c} in \cite{A} and \cite{A2}
using the results obtained Adamovi\'{c} and Milas 
\cite{AM}, Feigin, Semikhatov and Tipunin \cite{FST} and 
Doerrzapf \cite{D}:

\begin{thm}\label{ad}
The vertex operator superalgebra $L_{\mathfrak{ns}(2)}(c, 0, 0)$ 
has finitely many irreducible 
modules and every module for $L_{\mathfrak{ns}(2)}(c, 0, 0)$ 
is completely reducible if
and only if 
$$
c=c_{m}=\frac{3m}{m+2}
$$
where $m$ is a nonnegative integer.
A set of representatives of the equivalence classes of irreducible modules
for $L_{\mathfrak{ns}(2)}(c_{m}, 0, 0)$ is 
$$
\{L_{\mathfrak{ns}(2)}(c_{m}, h_{m}^{j,k}, q_{m}^{j,k})\}_{j, k\in 
\mathbb{N}_{\frac{1}{2}}, \;0\le j, k, j+k<m}
$$
where $\mathbb{N}_{\frac{1}{2}}
=\{\frac{1}{2}, \frac{3}{2}, \frac{5}{2},\dots \}$
and
\begin{eqnarray*}
h_{m}^{j, k}&=&\frac{jk-\frac{1}{4}}{m+2},\nn
q_{m}^{j, k}&=&\frac{j-k}{m+2}
 \hspace{\fill}
\epfe
\end{eqnarray*}
\end{thm}

For any $m\ge 0$,
we call the vertex operator algebra $L_{\mathfrak{ns}(2)}(c_{m}, 0, 0)$ a
{\it unitary minimal $N=2$ superconformal vertex operator superalgebra}.

\begin{prop} \label{fusion}
Let $j_{i}, k_{i}\in \mathbb{N}_{\frac{1}{2}}$, 
$i=1, 2, 3$, satisfying $0\le j_{i}, k_{i}, j_{i}+k_{i}<m$ 
and
$\mathcal{ Y}$ an intertwining operator of type
\begin{equation}\label{type}
{L_{\mathfrak{n}\mathfrak{s}(2)}(c_{m}, h_{m}^{j_{3}, k_{3}}, 
q_{m}^{j_{3}, k_{3}})\choose 
L_{\mathfrak{n}\mathfrak{s}(2)}(c_{m}, h_{m}^{j_{1}, k_{1}}, 
q_{m}^{j_{1}, k_{1}})
L_{\mathfrak{n}\mathfrak{s}(2)}(c_{m}, h_{m}^{j_{2}, k_{2}}, 
q_{m}^{j_{2}, k_{2}})}.
\end{equation} 
Then we have:

\begin{enumerate}
\item\label{2.2-1}  For any 
$w_{(1)}\in L_{\mathfrak{n}\mathfrak{s}(2)}(c_{m}, h_{m}^{j_{1}, k_{1}}, 
q_{m}^{j_{1}, k_{1}})$ 
and $w_{(2)}\in L_{\mathfrak{n}\mathfrak{s}(2)}(c_{m}, h_{m}^{j_{2}, k_{2}}, 
q_{m}^{j_{2}, k_{2}})$, 
$$\mathcal{Y}(w_{(1)}, x)w_{(2)}\in x^{h_{m}^{j_{3}, k_{3}}-h_{m}^{j_{1}, 
k_{1}}-h_{m}^{j_{2}, k_{2}}}L_{\mathfrak{n}\mathfrak{s}(2)}(c_{m}, 
h_{m}^{j_{3}, k_{3}}, 
q_{m}^{j_{3}, k_{3}})((x^{1/2})).$$

\item\label{2.2-2} Let $\Delta=h_{m}^{j_{3}, k_{3}}-h_{m}^{j_{1}, k_{1}}-
h_{m}^{j_{2}, k_{2}}$ and $w_{(i)}=\mathbf{1}_{c_{m}, h_{m}^{j_{i}, k_{i}}, 
q_{m}^{j_{i}, k_{i}}}$, $i=1, 2, 3$,
the lowest weight vectors
in $L_{\mathfrak{n}\mathfrak{s}(2)}(c_{m}, h_{m}^{j_{i}, k_{i}}, 
q_{m}^{j_{i}, k_{i}})$.
Then the map $\mathcal{Y}$ is uniquely 
determined by the maps
\begin{eqnarray*}
&(w_{(1)})_{-\Delta-1},&\\ 
&(G^{+}(-1/2)w_{(1)})_{-\Delta-1/2},&\\
&(G^{-}(-1/2)w_{(1)})_{-\Delta-1/2},&\\
&(G^{+}(-1/2)G^{-}(-1/2)w_{(1)})_{-\Delta}&
\end{eqnarray*}
{}from the one-dimensional subspace of $W_{2}$ spanned by 
$w_{(2)}$ to the one-dimensional subspace of $W_{3}$ spanned by 
$w_{(3)}$. That is,
if these maps are $0$,
then
$\mathcal{Y}=0$.

\item\label{2.2-3} If $q_{m}^{j_{3}, k_{3}}$ is not equal to one of the 
numbers $q_{m}^{j_{1}, k_{1}}+q_{m}^{j_{2}, k_{2}}$,
$q_{m}^{j_{1}, k_{1}}+q_{m}^{j_{2}, k_{2}}-1$ and 
$q_{m}^{j_{1}, k_{1}}+q_{m}^{j_{2}, k_{2}}+1$,
then the space 
$$\mathcal{V}^{L_{\mathfrak{n}\mathfrak{s}(2)}
(c_{m}, h_{m}^{j_{3}, k_{3}}, q_{m}^{j_{3}, k_{3}})}
_{L_{\mathfrak{n}\mathfrak{s}(2)}(c_{m}, 
h_{m}^{j_{1}, k_{1}}, q_{m}^{j_{1}, k_{1}})
L_{\mathfrak{n}\mathfrak{s}(2)}(c_{m}, 
h_{m}^{j_{2}, k_{2}}, q_{m}^{j_{2}, k_{2}})}$$ 
of intertwining operators of type (\ref{type})
is $0$. If $q_{m}^{j_{3}, k_{3}}= q_{m}^{j_{1}, k_{1}}
+q_{m}^{j_{2}, k_{2}}\pm 1$,
it is at most $1$-dimensional.
If $q_{m}^{j_{3}, k_{3}}=q_{m}^{j_{1}, k_{1}}+q_{m}^{j_{2}, k_{2}}$,
it
is at most $2$-dimensional.

\end{enumerate}
\end{prop}
\pf
Conclusion \ref{2.2-1} is clear since the three modules are 
irreducible.

Conclusion \ref{2.2-2} can be proved similarly to the 
proof of the similar statement in the $N=1$ case in \cite{HM}.
Here we give a different proof.

Suppose that  
\begin{eqnarray}
&(w_{(1)})_{-\Delta-1}w_{(2)},&\label{v1}\\
&(G^{+}(-1/2)w_{(1)})_{-\Delta-1/2}w_{(2)},&\label{v2}\\
&(G^{-}(-1/2)w_{(1)})_{-\Delta-1/2}w_{(2)},&\label{v3}\\
&(G^{+}(-1/2)G^{-}(-1/2)w_{(1)})_{-\Delta}w_{(2)}&\label{v4}
\end{eqnarray}
are all equal to $0$
but $\mathcal{Y}\neq 0$.

Using the associator formula (obtained by taking residue in 
$x_{1}$ in the Jacobi identity defining intertwining operators)
\begin{eqnarray*}
\mathcal{Y}(Y(u, x_{0})w, x_{2})&=&Y(u, x_{0}+x_{2})\mathcal{Y}(w, x_{2})\nn
&&-(-1)^{|u||w|}
\res_{x_{1}}x_{0}^{-1}\delta\left(\frac{x_{2}-x_{1}}{-x_{0}}\right)
\mathcal{Y}(w, x_{2})Y(u, x_{1})
\end{eqnarray*}
repeatedly, we see that $\mathcal{Y}= 0$
if $\mathcal{Y}(w_{(1)}, x)=0$. Thus $\mathcal{Y}(w_{(1)}, x)\neq 0$.

Using the commutator formula (obtained by taking residue in 
$x_{0}$ in the Jacobi identity defining intertwining operators)
\begin{eqnarray*}
\lefteqn{Y(u, x_{1})\mathcal{Y}(w, x_{2})-(-1)^{|u||w|}
\mathcal{Y}(w, x_{2})Y(w, x_{1})}\nn
&&=\res_{x_{0}}x_{2}^{-1}\delta\left(\frac{x_{1}-x_{0}}{x_{2}}\right)
\mathcal{Y}(Y(u, x_{0})w, x_{2})
\end{eqnarray*}
together with the $N=2$ Neveu-Schwarz algebra relations, the 
$L(-1)$-derivative property and the 
definition of the lowest weight vector $w_{(1)}$ repeatedly, 
we see that $\mathcal{Y}(w_{(1)}, x)=0$ if (\ref{v1})--(\ref{v4}) 
are all equal to $0$.
Thus these four vectors cannot be all $0$. We have a contradiction.

We prove Conclusion \ref{2.2-3} now.  
By Conclusion \ref{2.2-2}, we need only estimate the number
of nonzero vectors in the set of four vectors
(\ref{v1})--(\ref{v4}). 

We need the following:

\begin{lemma}\label{u1c}
The following equality hold:
\begin{equation} 
q_{m}^{j_{3}, k_{3}}(w_{(1)})_{-\Delta-1}w_{(2)}=
(q_{m}^{j_{1}, k_{1}}+q_{m}^{j_{2}, k_{2}})
(w_{(1)})_{-\Delta-1}w_{(2)}, \label{u1c1}
\end{equation}
\begin{eqnarray}
\lefteqn{q_{m}^{j_{3}, k_{3}}(G^+(-1/2)w_{(1)})_{-\Delta-1/2}w_{(2)}}\nn
&&=(q_{m}^{j_{1}, k_{1}}+q_{m}^{j_{2}, k_{2}}+1)
(G^+(-1/2)w_{(1)})_{-\Delta-1/2}w_{(2)},\label{u1c2}
\end{eqnarray}
\begin{eqnarray}
\lefteqn{q_{m}^{j_{3}, k_{3}}(G^+(-1/2)w_{(1)})_{-\Delta-1/2}w_{(2)}}\nn
&&= (q_{m}^{j_{1}, k_{1}}+q_{m}^{j_{2}, k_{2}}-1)
(G^+(-1/2)w_{(1)})_{-\Delta-1/2}w_{(2)},\label{u1c3}
\end{eqnarray}
\begin{eqnarray}\label{u1c4}
\lefteqn{q_{m}^{j_{3}, k_{3}}(G^+(-1/2)G^-(-1/2)
w_{(1)})_{-\Delta}w_{(2)}}\nn
&&=(q_{m}^{j_{1}, k_{1}}+q_{m}^{j_{2}, k_{2}})
(G^+(-1/2)G^-(-1/2)w_{(1)})_{-\Delta}w_{(2)}\nn
&&\quad +(2h_{m}^{j_{1}, k_{1}}+q_{m}^{j_{1}, k_{1}})
(w_{(1)})_{-\Delta-1}w_{(2)}.
\end{eqnarray}
\end{lemma}
\pf
A straightforward calculation gives
$$J(0)(w_{(1)})_{-\Delta-1}w_{(2)}=(q_{m}^{j_{1}, k_{1}}+q_{m}^{j_{2}, k_{2}})
(w_{(1)})_{-\Delta-1}w_{(2)},$$
\begin{eqnarray*}
\lefteqn{J(0)(G^+(-1/2)w_{(1)})_{-\Delta-1/2}w_{(2)}}\nn
&&=(q_{m}^{j_{1}, k_{1}}+q_{m}^{j_{2}, k_{2}}+1)
(G^+(-1/2)w_{(1)})_{-\Delta-1/2}w_{(2)},
\end{eqnarray*}
\begin{eqnarray*}
\lefteqn{J(0)(G^+(-1/2)w_{(1)})_{-\Delta-1/2}w_{(2)}}\nn
&&= (q_{m}^{j_{1}, k_{1}}+q_{m}^{j_{2}, k_{2}}-1)
(G^+(-1/2)w_{(1)})_{-\Delta-1/2}w_{(2)},
\end{eqnarray*}
\begin{eqnarray*}
\lefteqn{J(0)(G^+(-1/2)G^-(-1/2)
w_{(1)})_{-\Delta}w_{(2)}}\nn
&&=(q_{m}^{j_{1}, k_{1}}+q_{m}^{j_{2}, k_{2}})
(G^+(-1/2)G^-(-1/2)w_{(1)})_{-\Delta}w_{(2)}\nn
&&\quad +(2h_{m}^{j_{1}, k_{1}}+q_{m}^{j_{1}, k_{1}})
(w_{(1)})_{-\Delta-1}w_{(2)}.
\end{eqnarray*}
But on the other hand, note that (\ref{v1})--(\ref{v4}) are all 
(zero or nonzero) constant multiple of $w_{(3)}$ and thus all have 
$U(1)$ charge $q_{m}^{j_{3}, k_{3}}$. 
So we have (\ref{u1c1})--(\ref{u1c4}).
\epfv

We prove Conclusion \ref{2.2-3} using this lemma now.
If $q_{m}^{j_{3}, k_{3}}$ is not equal to one of the 
numbers $q_{m}^{j_{1}, k_{1}}+q_{m}^{j_{2}, k_{2}}$
$q_{m}^{j_{1}, k_{1}}+q_{m}^{j_{2}, k_{2}}-1$ and 
$q_{m}^{j_{1}, k_{1}}+q_{m}^{j_{2}, k_{2}}+1$,
then {}from (\ref{u1c1})--(\ref{u1c4}), we conclude 
that (\ref{v1})--(\ref{v4}) are all equal to $0$.
Thus the space of intertwining operators 
is $0$. 

If $q_{m}^{j_{3}, k_{3}}=q_{m}^{j_{1}, k_{1}}+q_{m}^{j_{2}, k_{2}}+1$,
then by  (\ref{u1c1}), (\ref{u1c3}) and (\ref{u1c4}), 
(\ref{v1}), (\ref{v3})
and (\ref{v4})
must be $0$. Thus there is at most one nonzero vector 
(\ref{v2}). So the 
dimension is at most $1$.

Similarly in the case of 
$q_{m}^{j_{3}, k_{3}}=q_{m}^{j_{1}, k_{1}}+q_{m}^{j_{2}, k_{2}}-1$,
we can show that $(w_{(1)})_{-\Delta-1}w_{(2)}$,
(\ref{v1}), (\ref{v2})
and (\ref{v4}) must be $0$
and thus the dimension is at most $1$.

If $q_{m}^{j_{3}, k_{3}}=q_{m}^{j_{1}, k_{1}}+q_{m}^{j_{2}, k_{2}}$,
then by (\ref{u1c2}) and (\ref{u1c3}), 
(\ref{v2}) and
(\ref{v3}) must both be $0$.
Thus we have at most two nonzero vecotrs (\ref{v1})
and (\ref{v4}) and the dimension
is at most $2$. 
\epf

\begin{rema}
{\rm In the proofs above we
do not use the particular properties, except the irreducibility,
of $L(c_{m}, h_{m}^{j_{i},
k_{i}},
q_{m}^{j_{i}, k_{i}})$, $i=1, 2, 3$.
Thus the conclusions of  Proposition \ref{fusion} 
remain true
if we replace $c_m$ by an arbitrary $c$ and
$L(c_m,h_{m}^{j_i,k_i},q_m^{j_i,k_i})$, $i=1,2,3$ by
$L(c,0,0)$-modules $L(c,h_i,q_i)$, $i=1,2,3$, if
$L(c,h_i,q_i)$ are irreducible.}
\end{rema}

\begin{defn}
{\rm An irreducible module for 
$L_{\mathfrak{n}\mathfrak{s}(2)}(c, 0, 0)$ is said to be chiral
(anti-chiral) if 
$$G^+(-1/2)w=0, \ \ (G^-(-1/2)w=0)$$
where $w$ is a nonzero lowest weight vector of the module.}
\end{defn}

Note that in the case $c=c_m$ we have only
finitely many
chiral (anti-chiral) modules.

\begin{cor}\label{fusion-ch}
Assume that $L_{\mathfrak{n}\mathfrak{s}(2)}(c_{m}, h_{m}^{j_{1}, k_{1}}, 
q_{m}^{j_{1}, k_{1}})$ is chiral or anti-chiral.
Then the dimension of the space 
$$\mathcal{V}^{L_{\mathfrak{n}\mathfrak{s}(2)}
(c_{m}, h_{m}^{j_{3}, k_{3}})}_{L_{\mathfrak{n}\mathfrak{s}(2)}(c_{m}, 
h_{m}^{j_{1}, k_{1}})L_{\mathfrak{n}\mathfrak{s}(2)}(c_{m}, 
h_{m}^{j_{2}, k_{2}})}$$
is at most $1$.
\end{cor}
\pf
Assume that $L_{\mathfrak{n}\mathfrak{s}(2)}(c_{m}, h_{m}^{j_{1}, k_{1}}, 
q_{m}^{j_{1}, k_{1}})$ is chiral. Then 
$$(G^{+}(-1/2)w_{(1)})_{-\Delta-1/2}w_{2}=0.$$
We claim that in this case
\begin{equation}\label{g+g-1}
(G^{+}(-1/2)G^{-}(-1/2)w_{(1)})_{-\Delta}w_{(2)}
=2(w_{(1)})_{-\Delta-1}w_{(2)}.
\end{equation}
In fact, the commutator formula for $G^{+}(-1/2)$ and $G^{-}(-1/2)$
gives 
$$
G^{+}(-1/2)G^{-}(-1/2)=G^{-}(-1/2)G^{+}(-1/2)
+2L(-1).
$$
Thus 
\begin{eqnarray}\label{g+g-2}
(G^{+}(-1/2)G^{-}(-1/2)w_{(1)})_{-\Delta}w_{2}&=&(G^{-}(-1/2)G^{+}(-1/2)
w_{(1)})_{-\Delta}w_{2}\nn
&&+2(L(-1)w_{(1)})_{-\Delta}w_{2}\nn
&=&2(L(-1)w_{(1)})_{-\Delta}w_{2}.
\end{eqnarray}
{}From the $L(-1)$-derivative property for intertwining operators, we obtain
$$(L(-1)w_{(1)})_{-\Delta}=(w_{(1)})_{-\Delta-1}.$$
Thus the right-hand side of (\ref{g+g-2}) becomes
$2(w_{(1)})_{-\Delta-1}w_{(2)}$, proving (\ref{g+g-1}).
On the other hand, {}from 
(\ref{u1c1}) and (\ref{u1c3}), we see that the vectors
$(w_{(1)})_{-\Delta-1}w_{2}$ and
$(G^{-}(-1/2)w_{(1)})_{-\Delta-1/2}w_{2}$ cannot be nonzero at the 
same time. Thus in this case, the dimension of the space spanned by
the four vectors (\ref{v1})--(\ref{v4}) is at most $1$.
Equivalently, the corollary is proved.
\epfv

\begin{rema}
{\rm Note that the operator $J(0)$ plays an essential role
in the proofs of Conclusion \ref{2.2-3} in Proposition \ref{fusion}
and Corollary \ref{fusion-ch}.}
\end{rema}

\begin{rema}
{\rm After the first version of the present paper was finished,
we received a preprint \cite{A3} from Adamovi\'{c} in which 
the fusion rules for $L_{\mathfrak{n}\mathfrak{s}(2)}(c_{m}, 0, 0)$
are calculated explicitly and a stronger complete reducibility theorem 
is proved. But for the purpose of the present paper, we shall not 
need these stronger results.}
\end{rema}

Combining Theorem \ref{ad} and the second or third conclusion of 
Proposition \ref{fusion}, we obtain:

\begin{cor}\label{rat}
The unitary minimal 
$N=2$ superconformal vertex operator superalgebras are rational in the 
sense of \cite{HL1}, that is, the following three conditions are 
satisfied: 

\begin{enumerate}

\item Every module for such an algebra is completely 
reducible.

\item There are only finitely many inequivalent irreducible modules for 
such an algebra.

\item The fusion rules among any three (irreducible) modules are 
finite.\epf

\end{enumerate}
\end{cor}

We also have:

\begin{prop}
Any finitely-generated 
lower truncated generalized 
module
$W$ for $L_{\mathfrak{n}\mathfrak{s}(2)}(c_{m},0, 0)$
 is an ordinary module.
\end{prop}
\pf
The proof is the same as the corresponding result in 
\cite{HM}. We repeat it here since it is simple. Suppose that $W$ is 
generated by a single vector $w \in W$. Then
by the Poincar\'{e}-Birkhoff-Witt theorem and 
the lower truncation condition, every
homogeneous  subspace of $U(\mathfrak{n}\mathfrak{s}(2))w$ is 
finite-dimensional, proving the result.
\epfv

Let $m_{i}$, $i=1, \dots, n$, be positive integers
and let 
$V=L_{\mathfrak{ns}(2)}(c_{m_{1}}, 0, 0)\otimes \cdots \otimes
L_{\mathfrak{ns}(2)}(c_{m_{n}}, 0, 0)$. 
{}From the trivial generalizations of the
results proved in \cite{FHL} and \cite{DMZ} to vertex operator
superalgebras, $V$ is a rational $N=2$ superconformal 
vertex operator superalgebra, a set
of representatives of equivalence classes of irreducible modules for
$V$ can be listed 
explicitly and the fusion rules for $V$ and 
can be calculated easily.

We introduce a class of $N=2$ superconformal vertex operator vertex
operator superalgebras:

\begin{defn}
{\rm Let $m_{i}$, $i=1, \dots, n$, be positive integers.
An $N=2$ superconformal  vertex
operator superalgebra $V$ is
said to be {\it in the class $\mathcal{C}_{m_{1};
\dots; m_{n}}$}
if $V$ has a vertex operator subalgebra 
isomorphic to $L_{\mathfrak{n}\mathfrak{s}(2)}(c_{m_{1}},
0, 0)\otimes \cdots \otimes 
L_{\mathfrak{n}\mathfrak{s}(2)}(c_{m_{n}}, 0, 0)$.}
\end{defn}

\begin{prop}\label{2-7}
Let $V$ be an $N=2$ superconformal vertex
operator superalgebra in the class $\mathcal{C}_{m_{1};
\dots; m_{n}}$. Then any finitely-generated 
lower truncated generalized $V$-module
$W$ is an ordinary module.
\end{prop}
\pf
The proof is similar to the proofs of Proposition 3.7 in \cite{H2}
and Proposition 2.7 in \cite{HM}.
Here we only point out the main difference. As in \cite{H2} and \cite{HM}, 
we discuss only the case $n=2$. Similar to 
the proofs of Proposition 3.7 in \cite{H2} and 
Proposition 2.7 in \cite{HM}, using the Jacobi identity,
the $N=2$ Neveu-Schwarz relations, in particular,
the formulas $[G^{+}(-1/2), G^{-}(-1/2)]=2L(-1)$, 
$(G^{+}(-1/2))^{2}=(G^{-}(-1/2))^{2}=0$, and Theorem 4.7.4 of \cite{FHL},
we can reduce our 
proof in the case of $n=2$ 
to the finite-dimensionality of the space spanned 
by the elements of the form 
\begin{eqnarray}\label{generators}
\lefteqn{A
(L(-1)^{l_{1}}G^{+}(-1/2)^{k_{1}}G^{-}(-1/2)^{k_{2}}
u_{(j)}^{(1)})_{j_{1}}
Bw_{(t)}^{(1)}}\nn
&&\quad\otimes C (L(-1)^{l_{2}}G^{+}(-1/2)^{k_{3}}G^{-}(-1/2)^{k_{4}}
u_{(j)}^{(2)})_{j_{2}}
Dw_{(t)}^{(2)},\nn
&&
\end{eqnarray}
$l_{1}, l_{2}\in \mathbb{N}$, $k_{1}, k_{2}, k_{3}, k_{4}=0, 1$,
$j_{1}, j_{2}\in \mathbb{Q}$, $t=1, \dots, c$,
$j=1, \dots, d$, where $A$ (or $C$) are 
products of operators of the forms $L(-a)_{1})$ (or $L(-a)_{2}$), 
$J(-a)_{1}$ (or $J(-a)_{2}$),
 $a\in \mathbb{Z}_{+}$, and
$G^{\pm}(-b)_{1}$ (or $G^{\pm}(-b)_{2}$), $b\in \mathbb{Z}/2$, $B$ (or $D$)
are products of operators of the forms $L(a)_{1}$ (or $L(a)_{2}$), 
$J(a)_{1}$ (or $J(a)_{2}$),
$a\in \mathbb{Z}_{+}$,
and $G^{\pm}(b)_{1}$ (or $G^{\pm}(b)_{2}$), $b\in \mathbb{Z}+\frac{1}{2}$,
where $u^{(i)}_{(j)}$, $j=1, \dots, d$, $i=1, 2$,
are elements of $V$ such that the 
$L_{\mathfrak{n}\mathfrak{s}(2)}(c_{m_{i}}, 0, 0)$-submodules 
generated by them are isomorphic to 
$L(c_{m_{i}}, h_{m_{i}}^{r_{j},s_{j}}, q_{m_{i}}^{r_{j},s_{j}})$ 
for some $r_{j},s_{j}\in 
\mathbb{N}_{\frac{1}{2}}$ satisfying $0\le r_{j},s_{j}, r_{j}+s_{j}<m_{i}$
with the images of 
$u^{(i)}_{(j)}$, $j=1, \dots, d$, $i=1, 2$, as the lowest weight vectors 
and such that $V$ is isomorphic to the direct sum of these submodules,
and where
$w_{(t)}^{(i)}$, $t=1, \dots, c$, $i=1, 2$, are homogeneous
elements of some irreducible 
$L_{\mathfrak{n}\mathfrak{s}(2)}(c_{m_{i}}, 0, 0)$-modules. 
The remaining argument in the proof is the same as the 
corresponding parts in \cite{H2} and \cite{HM} and is omitted here.
\epfv

\renewcommand{\theequation}{\thesection.\arabic{equation}}
\renewcommand{\thethm}{\thesection.\arabic{thm}}
\setcounter{equation}{0}
\setcounter{thm}{0}

\section{The convergence and the extension property
for products of intertwining operators}

In this section, we study products of intertwining operators for the
unitary minimal $N=2$ superconformal vertex operator superalgebras. The 
main result is the folllwing:

\begin{thm}\label{cep2}
Let $m$ be a positive integer. Then
intertwining operators for the vertex operator 
superalgebra $L(c_m,0,0)$ satisfy the convergence
and extension property for products of intertwining operators
introduced in \cite{H1}.
\end{thm}

An immediate consequence is a similar result for the vertex operator
algebras in the class $\mathcal{C}_{m_{1}; \dots; m_{n}}$. See Theorem
\ref{cep1}. 

Instead of proving Theorem \ref{cep2} by
deriving differential equations with regular singularities satisfied
by the matrix elemenets of products of intertwining operators of
lowest weight vectors, as in \cite{H1}, \cite{HL5.5} and \cite{HM}, we
use the so-called anti-Kazama-Suzuki mapping \cite{FST}, which reduces
the problem to the study of intertwining operators for a vertex
operator algebra constructed {}from an
$\widehat{\mathfrak{s}\mathfrak{l}}_2$ integrable lowest weight
representation.

First we give some auxiliary constructions 
and discuss some results obtained using the so-called
``anti-Kazama-Suzuki mapping''
(introduced in \cite{FST}).

Fix a  positive integer $m$. 
As before, we let $c_m=\frac{3m}{m+2}$
and denote by $L_{\mathfrak{n}\mathfrak{s}(2)}(c_m,0,0)$ 
the  unitary minimal $N=2$ superconformal vertex operator superalgebra.

Let $L$ be a rank one lattice generated by $\alpha$ with the
bilinear form $\langle \cdot, \cdot\rangle$ given by 
$$\langle \alpha,\alpha\rangle=-1,$$
and let  $l=L \otimes_{\mathbb{Z}} \mathbb{C}$.
As in \cite{FLM}, we have a vertex superalgebra
$$V_L \cong S(\hat{l}_-) \otimes \mathbb{C}[L].$$
Note that $V_{L}$ is super since $L$ is odd and $V_{L}$ does not 
satisfy the grading-restriction conditions for the grading obtained
{}from the usual Virasoro element for lattice vertex algebras since 
the bilinear form is not positive definite. 

Let 
$$\tilde{\omega}_{V_L}=-\frac{\alpha(-1)^2}{2}+\frac{\alpha(-2)}{2}.$$
It can be verified easily that the component operators of the vertex
operator associated to $\tilde{\omega}_{V_L}$ satisfy the Virasoro
relations with central charge $4$. In particular, the component
operator for the $-2$-nd power of $x$ gives a $\mathbb{Z}/2$-grading for
$V_{L}$. With this grading, $V_{L}$ is a $\mathbb{Z}/2$-graded vertex
superalgebra. However, it is easy to see that
this grading is not truncated {}from below. Thus
$V_{L}$ with $\tilde{\omega}_{V_L}$ as its Virasoro element fails to
be a vertex operator superalgebra.

We also need the following construction of the so-called ``Liouville
scalar model:'' Let $\hat{\mathbb{R}}$ be the Lie algebra with a basis
$a(n)$, $n \in \mathbb{Z}$, and $d$, satisfying the bracket relations
\begin{eqnarray*} {[a(m),a(n)]}&=&m\delta_{m+n,0}d,\\ 
{[a(m), d]}&=&0
\end{eqnarray*} 
for $m, n\in \mathbb{Z}$.  Let $M(1,s)$ be the
corresponding irreducible highest weight module with central charge $1$
and highest weight $s$.  It is well-known that
$M(1,0)=S(\hat{\mathbb{R}})$ has a vertex operator algebra structure
with the Virasoro element $\omega_{M(1, 0)}=\frac{a(-1)^2}{2}$.
Consider the vertex algebra structure on $M(1,0)$ together with a
different Virasoro element
$$\tilde{\omega}_{M(1,0)}=\frac{a(-1)^2}{2}+\frac{ia(-2)}{2}.$$ 
A straightforward calculation shows that the vertex algebra
structure on $M(1,0)$ together with the Virasoro element
$\tilde{\omega}_{M(1,0)}$ is a vertex operator algebra with the 
central charge $4$ (see, for example, \cite{L} and \cite{FF}). 
We shall denote this vertex operator
algebra by $V_{\rm Liou}$ (the vertex operator algebra associated to the
Liouville scalar model).  In addition, every $M(1,s)$, $s \in
\mathbb{C}$, is an irreducible $V_{\rm Liou}$-module and any $V_{\rm
Liou}$-module on which $a(0)$ acts semisimply is completely reducible.
We will work only with such modules $V_{\rm Liou}$-module, which are
enough for our purposes. 

The anti-Kazama-Suzuki mapping gives us a structure of an
$\widehat{\mathfrak{sl}}_2$-module on 
$L_{\mathfrak{ns}(2)}(c_m,h_{m}^{j, k},q_{m}^{j, k}) \otimes V_L$
for $j, k\in \mathbb{N}_{\frac{1}{2}}$, $0<j, k, j+k<m$.  
Consider the vectors (as
in \cite{FST} and \cite{A})
\begin{eqnarray*}
\mathfrak{e}&=&G^{+}(-3/2)\mathbf{1}_{c_{m}, 0, 0}\otimes e^{-\alpha},\\
 \mathfrak{f}&=&\frac{m+2}{2} G^{-}(-3/2)\mathbf{1}_{c_{m}, 0, 0}
 \otimes e^{\alpha},\\
 \mathfrak{h}&=&-m\mathbf{1}_{c_{m}, 0, 0}\otimes
\alpha(-1)+(m+2)J(-1)\mathbf{1}_{c_{m}, 0, 0}\otimes e^{0}
\end{eqnarray*}
in $L_{\mathfrak{ns}(2)}(c_m, 0, 0) \otimes V_L$.
Then the vertex operators
$Y(\mathfrak{e},x)$, $Y(\mathfrak{f},x)$ and $Y(\mathfrak{h},x)$ for
the $L_{\mathfrak{ns}(2)}(c_m,0,0) \otimes V_L$-module 
$L_{\mathfrak{ns}(2)}(c_m,h_{m}^{j, k},q_{m}^{j, k}) \otimes V_L$ 
give
a representation of $\widehat{\mathfrak{sl}}_2$ of
 level $m$ on $L_{\mathfrak{ns}(2)}(c_m,h_{m}^{j, k},q_{m}^{j, k}) 
\otimes V_L$. The main
observation in \cite{FST} is that 
$L_{\mathfrak{ns}(2)}(c_m,h_{m}^{j, k},q_{m}^{j, k}) \otimes
V_L$ is completely reducible as an $\widehat{\mathfrak{sl}}_2$-module.
In the special case $h_{m}^{j, k}=q_{m}^{j, k}=0$, 
we obtain a vertex subalgebra of
$L_{\mathfrak{ns}(2)}(c_m,0,0) \otimes V_L$ which is isomorphic 
as a vertex algebra to the
underlying vertex algebra of the vertex operator algebra
$L_{\widehat{\mathfrak{sl}}(2)}(m, 0)$ on the integrable highest-weight
$\mathfrak{sl}(2)$-module of level $m$ and highest weight $0$ (as in
\cite{A}).  The Virasoro elelement for $L_{\widehat{\mathfrak{sl}}(2)}(m,
0)$ is given by the Sugawara-Segal construction, and if we identify 
this vertex subalgebra with
$L_{\widehat{\mathfrak{sl}}(2)}(m, 0)$, then
\begin{eqnarray*}
\omega_{L_{\widehat{\mathfrak{sl}}(2)}(m,0)}
&=&\omega_{L_{\mathfrak{ns}(2)}(c_m,0,0)}\otimes e^{0}
+\frac{m+2}{4}J(-1)^2\mathbf{1}_{c_{m}, 0, 0}\otimes e^{0}\nn
&&-
\frac{m+2}{2}J(-1)\mathbf{1}_{c_{m}, 0, 0}\otimes \alpha(-1)
+ \frac{m}{4}\mathbf{1}_{c_{m}, 0, 0}
\otimes \alpha(-1)^2,
\end{eqnarray*}
where $\omega_{L_{\widehat{\mathfrak{sl}}(2)}(m,0)}$ and 
$\omega_{L_{\mathfrak{ns}(2)}(c_m,0,0)}$ are the Virasoro elements for 
$L_{\widehat{\mathfrak{sl}}(2)}(m,0)$ and $L_{\mathfrak{ns}(2)}(c_m,0,0)$,
respectively.
We see that under the isomorphism {}from this vertex subalgebra of
$L_{\mathfrak{ns}(2)}(c_m,0,0) \otimes V_L$ to 
$L_{\widehat{\mathfrak{sl}}(2)}(m, 0)$, the
Virasoro element in $L_{\widehat{\mathfrak{sl}}(2)}(m, 0)$ is not the
image of the Virasoro element of $L_{\mathfrak{ns}(2)}(c_m,0,0) \otimes V_L$. 

To get the correct Virasoro element (as in \cite{FST}),
we consider the vertex subalgebra of 
$L_{\mathfrak{ns}(2)}(c_m,0,0) \otimes V_L$
generated by the element
$$\rho=\sqrt{\frac{m+2}{2}}(J(-1)\mathbf{1}_{c_{m}, 0, 0}\otimes e^{0}-
\mathbf{1}_{c_{m}, 0, 0}\otimes \alpha(-1)).$$
It is straightforward to verify that this vertex subalgebra is 
actually isomorphic to $V_{\rm Liou}$. Straightforward calculations
also show that $Y(\rho,x)$ commutes
with $\widehat{\mathfrak{sl}}_2$ generators. So
$L_{\widehat{\mathfrak{sl}}(2)}(m,0) \otimes V_{\rm Liou}$ 
is isomorphic to a vertex subalgebra 
of $L_{\mathfrak{ns}(2)}(c_m,0,0) \otimes V_L$ as well and 
we shall, for convenience, identify 
$L_{\widehat{\mathfrak{sl}}(2)}(m,0) \otimes V_{\rm Liou}$ with this  
vertex subalgebra.
It is easy to see that the Virasoro element of 
$L_{\widehat{\mathfrak{sl}}(2)}(m,0) \otimes V_{\rm Liou}$ 
is identified with the Virasoro element
of $L_{\mathfrak{ns}(2)}(c_m,0,0)\otimes V_L$.
Thus $L_{\widehat{\mathfrak{sl}}(2)}(m,0) \otimes V_{\rm Liou}$ 
is a vertex operator subalgebra.

Now we recall a key result {}from \cite{FST} and \cite{FSST}
(in the case of unitary modules), 
slightly reformulated in the language of vertex operator
algebras:
\begin{thm} \label{fst}
As a generalized $L_{\widehat{\mathfrak{sl}}(2)}(m,0) \otimes V_{\rm Liou}$-module, 
$L_{\mathfrak{ns}(2)}(c_m,h,q) \otimes V_L$ decomposes
as
$$\oplus_{k \in \{0,1,\ldots,m \}} 
\oplus_{s \in I_s}  L_{\widehat{\mathfrak{sl}}(2)}(m,k) \otimes
M(1,s).$$
where $s$ runs through certain infinite index set $I_s$.
\end{thm}

We know that $L_{\widehat{\mathfrak{sl}}(2)}(m,0)$
is rational (see \cite{FZ}). Although $V_{\rm Liou}$ is not 
rational,
any module on which $a(0)$ acts semisimply is completely 
reducible, as we mentioned above. 
Any irreducible module for $L_{\widehat{\mathfrak{sl}}(2)}(m,0)$ is 
isomorphic to $L_{\widehat{\mathfrak{sl}}(2)}(m,i)$ for some 
$i\in \{1, \dots, m\}$ and any
irreducible module for $V_{\rm Liou}$ is isomorphic to 
$M(1,s)$ for some $s\in \mathbb{C}$. Thus by the result 
in \cite{FHL}
on modules for a tensor product of  vertex operator algebras,
any irreducible module for $L_{\widehat{\mathfrak{sl}}(2)}(m,0)\otimes V_{\rm Liou}$ is isomorphic to
$L_{\widehat{\mathfrak{sl}}(2)}(m,i) \otimes M(1,s)$ for some $i \in \{0,...,m\}$,
$s \in {\mathbb C}$. 

Proposition 2.7  in \cite{DMZ} and its proof can be generalized trivially
to the case that one of the vertex operator algebra is
an irrational vertex operator algebra like $V_{\rm
Liou}$, such that in particular, 
any $L_{\widehat{\mathfrak{sl}}(2)}(m,0) \otimes V_{\rm
Liou}$-module with
$\mathbf{1}_{L_{\widehat{\mathfrak{sl}}(2)}(m,0)} \otimes a(0)$ 
($\mathbf{1}_{L_{\widehat{\mathfrak{sl}}(2)}(m,0)}$ being the 
vacuum vector of $L_{\widehat{\mathfrak{sl}}(2)}(m,0)$)
acting semi-simply is completely reducible.  Now suppose that $M$ is 
such an $L_{\widehat{\mathfrak{sl}}(2)}(m,0) \otimes V_{\rm Liou}$-module.  
Then it follows that
$M$ is completely reducible. So it has a decomposition 
\begin{equation}
\label{decomposition} \oplus_{\beta \in B} M_\beta
\end{equation}
where $B$ is an index set. 
Note
that here the sum might be infinite (comparing with the rational case
where this sum is always finite).  Since any irreducible 
$L_{\widehat{\mathfrak{sl}}(2)}(m,0) \otimes
V_{\rm Liou}$-module is isomorphic to 
$L_{\widehat{\mathfrak{sl}}(2)}(m,i) \otimes M(1,s)$ for some
$i \in \{0,...,m\}$, $s \in {\mathbb C}$, $M_{\beta}$ for any $\beta\in B$
is isomorphic to such a module.

We need the following:
\begin{lemma} \label{facint}
Let $\mathcal{Y}$ be an intertwining operator of 
type 
$${L_{\widehat{\mathfrak{sl}}(2)}(m,i_3) \otimes M(1,s_3) \choose
L_{\widehat{\mathfrak{sl}}(2)}(m,i_1) \otimes M(1,s_1) \ \ L_{\widehat{\mathfrak{sl}}(2)}(m,i_2) \otimes M(1,s_2)}.$$
Then  
$$\mathcal{Y}=\mathcal{Y}' \otimes \mathcal{Y}''$$
where $\mathcal{Y}'$ and $\mathcal{Y}''$ are intertwining
operators of types ${L_{\widehat{\mathfrak{sl}}(2)}(m,i_3) \choose 
L_{\widehat{\mathfrak{sl}}(2)}(m,i_1)
L_{\widehat{\mathfrak{sl}}(2)}(m,i_2)}$ and
${M(1,s_3) \choose M(1,s_1)  M(1,s_2)}$, respectively.
In particular, all fusion rules for irreducible modules
for $L_{\widehat{\mathfrak{sl}}(2)}(m,0) \otimes
V_{\rm Liou}$ are $0$ or $1$.
\end{lemma}
\pf
It is enough to show that there is a linear injective map
{}from 
$$\mathcal{V}^{L_{\widehat{\mathfrak{sl}}(2)}(m,i_3) 
\otimes M(1,s_3)}_{(L_{\widehat{\mathfrak{sl}}(2)}(m,i_1) \otimes M(1,s_1))
(L_{\widehat{\mathfrak{sl}}(2)}(m,i_2) \otimes M(1,s_2))},$$
the space of intertwining operators
of type 
$${L_{\widehat{\mathfrak{sl}}(2)}(m,i_3) \otimes M(1,s_3)
\choose L_{\widehat{\mathfrak{sl}}(2)}(m,i_1) \otimes M(1,s_1)\ \ 
L_{\widehat{\mathfrak{sl}}(2)}(m,i_2) \otimes M(1,s_2)},$$
to
$$\mathcal{V}^{L_{\widehat{\mathfrak{sl}}(2)}(m,i_3)}
_{ L_{\widehat{\mathfrak{sl}}(2)}(m,i_1)
L_{\widehat{\mathfrak{sl}}(2)}(m,i_2)}
 \otimes \mathcal{V}^{M(1,s_3)}_{M(1,s_1) M(1,s_2)},$$ where
$\mathcal{V}^{L_{\widehat{\mathfrak{sl}}(2)}(m,i_3)}
_{ L_{\widehat{\mathfrak{sl}}(2)}(m,i_1)
L_{\widehat{\mathfrak{sl}}(2)}(m,i_2)}$ and
$\mathcal{V}^{M(1,s_3)}_{M(1,s_1) M(1,s_2)}$ are the space of
intertwining operators of type 
${L_{\widehat{\mathfrak{sl}}(2)}(m,i_3)\choose
L_{\widehat{\mathfrak{sl}}(2)}(m,i_1)
L_{\widehat{\mathfrak{sl}}(2)}(m,i_2)}$ and
${M(1,s_3)\choose M(1,s_1) M(1,s_2)}$, respectively.  
But this follows {}from Proposition 2.10 in \cite{DMZ} which in turn 
is a consequence of a result in \cite{FHL} on irreducible modules for 
a tensor product vertex operator algebra and a result 
in \cite{FZ} giving an isomorphism between a space of 
intertwining operators and a certain vector space.
Since
the fusion rules for irreducible 
$L_{\widehat{\mathfrak{sl}}(2)}(m,0)$-modules
are $0$ and $1$
(see \cite{FZ}) and the same is true for irreducible $V_{\rm Liou}$-modules,
the lemma is proved. \epfv

\begin{rema}
{\rm Note that  this result should be very usefull in calculating
the fusion algbera for $L_{\mathfrak{ns}(2)}(c_m,0,0)$ since every
such intertwining operator for $L_{\mathfrak{ns}(2)}(c_m,0,0)$ factors
as a tensor product of an intertwining operator
 for vertex operator algebra $L_{\widehat{\mathfrak{sl}}(2)}(m,0)$ 
and an intertwining operator for the vertex operator algebra
associated to the Heisenberg
algebra. In the present paper, the exact values of the fusion rules
are not what we are interested and thus we shall not calculate
them here. We leave this calculation to the
interested readers. (After the first version of the present paper was
finished, we received a preprint \cite{A3} from Adamovi\'{c} in which 
the fusion rules for $L_{\mathfrak{ns}(2)}(c_m,0,0)$ are calculated
explicitly.)}
\end{rema}
 
Using Lemmas \ref{facint}, we obtain:
\begin{prop} \label{propfinite}
For fixed $i_{1}, i_{2}\in \{1, \dots, m\}$
and $s_{1}, s_{2}\in \mathbb{C}$,
if $s_{3}\ne s_{1}+s_{2}$, the space 
$$\mathcal{V}_{(L_{\widehat{\mathfrak{sl}}(2)}(m,i_1) 
\otimes M(1,s_1))(L_{\widehat{\mathfrak{sl}}(2)}(m,i_{3}) 
\otimes M(1,s_2))}^{L_{\widehat{\mathfrak{sl}}(2)}(m,i_{3}) 
\otimes M(1,s_3)}$$
is $0$. In particular,
there are only finitely many pairs $(i_{3}, s_{3})
\in \{1, \dots,m\}\times \mathbb{C}$ such that 
the space 
$$\mathcal{V}_{(L_{\widehat{\mathfrak{sl}}(2)}(m,i_1) \otimes 
M(1,s_1))(L_{\widehat{\mathfrak{sl}}(2)}(m,i_{3}) 
\otimes M(1,s_2))}^{L_{\widehat{\mathfrak{sl}}(2)}(m,i_{3}) 
\otimes M(1,s_3)}$$
are not $0$. 
\end{prop}
\pf
Let 
$$\mathcal{Y}\in \mathcal{V}_{(L_{\widehat{\mathfrak{sl}}(2)}
(m,i_1) \otimes M(1,s_1))(L_{\widehat{\mathfrak{sl}}(2)}(m,i_{3}) 
\otimes M(1,s_2))}^{L_{\widehat{\mathfrak{sl}}(2)}(m,i_{3}) 
\otimes M(1,s_3)}.$$
By Lemma \ref{facint}, 
$$\mathcal{Y}=\mathcal{Y}'\otimes \mathcal{Y}''$$
where $\mathcal{Y}'$ and $\mathcal{Y}''$ are intertwining
operators of types ${L_{\widehat{\mathfrak{sl}}(2)}(m,i_3) 
\choose L_{\widehat{\mathfrak{sl}}(2)}(m,i_1) 
L_{\widehat{\mathfrak{sl}}(2)}(m,i_2)}$ and
${M(1,s_3) \choose M(1,s_1)  M(1,s_2)}$, respectively.
It is clear that if $s_{3}\ne s_{1}+s_{2}$, 
$\mathcal{Y}''=0$, proving the result.
\epfv

{\it Proof of Theorem \ref{cep2}}:\hspace{2ex}
Let $W_l$, $l=1,...,5$, be irreducible 
$L_{\mathfrak{ns}(2)}(c_m,0,0)$-modules and
 $\mathcal{Y}_1$
and $\mathcal{Y}_2$ intertwining operators of types
${W_4 \choose W_1 W_5}$ and ${W_5 \choose W_2 W_3}$, respectively.
Consider the (formal) matrix coefficients
\begin{equation} \label{matrix1}
\langle w'_{(4)}, \mathcal{Y}_1(w_{(1)},x_1)\mathcal{Y}_2(w_{(2)},x_2)
w_{(3)} \rangle,
\end{equation}
where $w_{(l)} \in W_l$, $l=1,2,3$ and $w'_{(4)}\in W'_{4}$.
We shall identify $W_{l}$, $l=1, 2, 3$, with $W_{l}\otimes e^{0}$ in
$W_{l} \otimes V_L$ and $W'_{4}$ with $W'_{4}\otimes e^{0}$
in $W'_{4} \otimes V_L$. In particular,
we use the same notations $w_{(l)}$, $l=1,2,3$, to denote  
$w_{(l)} \otimes e^{0}$,
and $w'_{(4)}$ to denote $w'_{(4)} \otimes e^{0}$.

We extend  intertwining operators $\mathcal{Y}_1$ and
$\mathcal{Y}_2$ uniquely to intertwining operators 
(denoted by the same notations $\mathcal{Y}_1$ and $\mathcal{Y}_2$)
of type
\begin{equation}
{W_4 \otimes V_L \choose W_1\otimes V_L \ \  W_5\otimes V_L},
\label{type1}
\end{equation}
and
\begin{equation}
{W_5 \otimes V_L \choose W_2\otimes V_L \ \ W_3\otimes V_L},
\label{type2}
\end{equation}
respectivly. 
By  Theorem \ref{fst}, $W_l \otimes V_L$, $l=1, 2, 3$, and 
$W'_{4}\otimes W'_{4}$ 
are
generalized modules for $L_{\widehat{\mathfrak{sl}}(2)}(m,0) 
\otimes V_{\rm Liou}$ and are 
completely reducible. 
So $w_{(l)}=\sum_{k=1}^{p_{l}}w_{(l)}^{(k)}$, $l=1, 2, 3$, and 
$w'_{(4)}=\sum_{k=1}^{p_{4}}w_{(4)}^{(k)}$ where 
$w_{(l)}^{(k)}$, $k=1, \dots, p_{i}$, $l=1, 2, 3, 4$, are 
elements of direct summands $M_{l}^{(k)}$ (irreducible 
$L_{\widehat{\mathfrak{sl}}(2)}(m,0) \otimes 
V_{\rm Liou}$-modules)
in $W_{i}$ for $i=1, 2, 3$ or $W'_{4}$ for $i=4$. 
Thus (\ref{matrix1}) is equal to 
\begin{equation}\label{matrix2}
\sum_{k_{1}=1}^{p_{1}}\sum_{k_{2}=1}^{p_{2}}\sum_{k_{2}=1}^{p_{3}}
\sum_{k_{4}=1}^{p_{4}}\langle w_{(4)}^{(k_{4})},
\mathcal{Y}^{k_{4}}_{k_{1}W_{5}}(w_{(1)}^{(k_{1})}, x_{1})
\mathcal{Y}^{W_{5}}_{k_{2}k_{3}}(w_{(2)}^{(k_{2})}, x_{2})
w_{(3)}^{(k_{3})}\rangle,
\end{equation}
where $\mathcal{Y}^{k_{4}}_{k_{1}W_{5}}$ and 
$\mathcal{Y}^{W_{5}}_{k_{2}k_{3}}$ are intertwining operators of 
types ${M_{4}^{(k_{4})}\choose M_{1}^{(k_{1})}W_{5}}$
and ${W_{5}\choose M_{2}^{(k_{2})}M_{3}^{(k_{3})}}$, respectively. 

By  Theorem \ref{fst}, $W_{5}$ is a completely reducible generalized 
$L_{\widehat{\mathfrak{sl}}(2)}(m,0) \otimes V_{\rm Liou}$-module. 
By Proposition \ref{propfinite}, 
(\ref{matrix2}) is equal to 
\begin{equation}\label{matrix3}
\sum_{k_{1}=1}^{p_{1}}\sum_{k_{2}=1}^{p_{2}}\sum_{k_{2}=1}^{p_{3}}
\sum_{k_{4}=1}^{p_{4}}\sum_{k_{5}=1}^{p_{5}}\langle w_{(4)}^{(k_{4})},
\mathcal{Y}^{k_{4}}_{k_{1}k_{5}}(w_{(1)}^{(k_{1})}, x_{1})
\mathcal{Y}^{k_{5}}_{k_{2}k_{3}}(w_{(2)}^{(k_{2})}, x_{2})
w_{(3)}^{(k_{3})}\rangle,
\end{equation}
where $\mathcal{Y}^{k_{4}}_{k_{1}k_{5}}$ and 
$\mathcal{Y}^{k_{5}}_{k_{2}k_{3}}$ are intertwining operators of 
types ${M_{4}^{(k_{4})}\choose M_{1}^{(k_{1})}M_{5}^{(k_{5})}}$
and ${M_{5}^{(k_{5})}\choose M_{2}^{(k_{2})}M_{3}^{(k_{3})}}$, respectively,
and $M_{5}^{(k_{5})}$, $k_{5}=1, \dots, p_{5}$, are irreducible 
irreducible $L_{\widehat{\mathfrak{sl}}(2)}(m,0) 
\otimes V_{\rm Liou}$-submodules of
$W_{5}$ as an $L_{\widehat{\mathfrak{sl}}(2)}(m,0) 
\otimes V_{\rm Liou}$-module.

By Proposition \ref{facint},
$$\mathcal{Y}^{k_{4}}_{k_{1}k_{5}}=(\mathcal{Y}^{k_{4}}_{k_{1}k_{5}})'
\otimes (\mathcal{Y}^{k_{4}}_{k_{1}k_{5}})''$$
and 
$$\mathcal{Y}^{k_{5}}_{k_{2}k_{3}}=(\mathcal{Y}^{k_{5}}_{k_{2}k_{3}})'
\otimes (\mathcal{Y}^{k_{5}}_{k_{2}k_{3}})'',$$
where $(\mathcal{Y}^{k_{4}}_{k_{1}k_{5}})'$ and 
$(\mathcal{Y}^{k_{5}}_{k_{2}k_{3}})'$ are intertwining operators for 
the vertex operator algebras
$L_{\widehat{\mathfrak{sl}}(2)}(m, 0)$ and $(\mathcal{Y}^{k_{4}}_{k_{1}k_{5}})''$ and 
$(\mathcal{Y}^{k_{5}}_{k_{2}k_{3}})''$ are intertwining operators for 
the vertex operator algebras $V_{\rm Liou}$. Thus (\ref{matrix3}) 
is equal to a finite sum of series of the form
\begin{eqnarray} \label{matrix}
\langle \tilde{w}_{(4)},
\widetilde{\mathcal{Y}}_{1}(\tilde{w}_{(1)}, x_{1})
\widetilde{\mathcal{Y}}_{2}(\tilde{w}_{(2)}, x_{2})
\tilde{w}_{(3)}\rangle\langle \tilde{w}_{(4)},
\widetilde{\widetilde{\mathcal{Y}}}_{1}(\tilde{\tilde{w}}_{(1)}, x_{1})
\widetilde{\widetilde{\mathcal{Y}}}_{2}(\tilde{\tilde{w}}_{(2)}, x_{2})
\tilde{\tilde{w}}_{(3)}\rangle,
\end{eqnarray}
where $\widetilde{\mathcal{Y}}_{1}$ and $\widetilde{\mathcal{Y}}_{2}$ 
are intertwining operators among irreducible modules
for $L_{\widehat{\mathfrak{sl}}(2)}(m,0)$ and 
$\widetilde{\widetilde{\mathcal{Y}}}_{1}$ and
$\widetilde{\widetilde{\mathcal{Y}}}_{2}$ are intertwining operators
among irreducible modules for $V_{\rm Liou}$. 

In \cite{HL5.5}, it was proved that intertwining operators for the
vertex operator algebra  $L_{\widehat{\mathfrak{sl}}(2)}(m,0)$ satisify 
the convergence and
extension property for products using  the Knizhnik-Zamolodchikov
eqautions. The convergence and
extension property for products of intertwining operators for the 
vertex operator algebra $V_{\rm Liou}$ can be proved trivially
by a straightford calculation. Using the 
convergence and extension properties for products of intertwining operators
for these vertex operator algebras and using
the fact proved above that (\ref{matrix1}) is a finite sum of 
series of the form (\ref{matrix}), we conclude that (\ref{matrix1})
is convergent when we substitute 
$e^{n_i \log z_i}$ for $x_{i}^{n}$ with $z_{1}, z_{2}\in \mathbb{C}$
satisfying $|z_1|>|z_2|>0$ and it can be analytically extended 
to an analytic function in the region
$|z_2|>|z_1-z_2|>0$ of the form 
\begin{equation} \label{extint}
\sum_{i=1}^j z_2^{r_i} (z_1-z_2)^{s_i}f_i \left(\frac{z_1-z_2}{z_2}\right).
\end{equation}

We still need to prove the following:
There exists $N$ (which does not depend on $w_{(1)}$ and $w_{(2)}$) 
such that
\begin{equation}\label{ineq}
{\rm wt}(w_{(1)})+{\rm wt}(w_{(2)})+s_i > N,
\end{equation}
for $i=1,...,j$.  The existence of $N$ follows (as in the $N=0$ and
$N=1$ cases in \cite{H2} and \cite{HM}, respectively) from an
induction argument for the $N=2$ superconformal algebra. Since any
$L_{\mathfrak{n}\mathfrak{s}(2)}(c_m,0,0)$-module is completely
reducible and since there are only finitely many irreducible
$L_{\mathfrak{n}\mathfrak{s}(2)}(c_m,0,0)$-modules, we need only prove
this existence in the case that $W_{1}$ and $W_{2}$ are irreducible.
When $w_{(1)}$ and $w_{(2)}$ are lowest weight vectors,
(\ref{matrix1}) is absolutely convergent in the region $|z_1|>|z_2|>0$
and can be analytically extended to an analytic function in in the
region $|z_2|>|z_1-z_2|>0$ of the form (\ref{extint}). We choose an
$N$ such that for these lowest weight vectors, (\ref{ineq}) holds. For
general $w_{(1})$ and $w_{(2)}$, we use induction instead of the proof
above to show that (\ref{matrix1}) in the region $|z_1|>|z_2|>0$
and can be analytically extended to an analytic function in in the
region $|z_2|>|z_1-z_2|>0$ of the form (\ref{extint}). In addition,
the induction also shows that (\ref{ineq}) holds for the $N$ we
choose.
\epfv

An immediate consequence of Theorem \ref{cep2}
is the following:

\begin{thm}\label{cep1}
Let $m_{i}$, $i=1, \dots, n$, be positive integers and 
$V$ an $N=2$ superconformal vertex operator superalgebra
in the class $\mathcal{C}_{m_{1};\dots, m_{n}}$. Then
intertwining operators for $V$ satisfy the convergence
and extension property for products of intertwining operators
introduced in \cite{H1}.\epf
\end{thm}

We omit the proof since it is the same as the corresponding 
result in \cite{H2} and \cite{HL5.5}.

\renewcommand{\theequation}{\thesection.\arabic{equation}}
\renewcommand{\thethm}{\thesection.\arabic{thm}}
\setcounter{equation}{0}
\setcounter{thm}{0}

\section{Intertwining operator superalgebras and 
vertex tensor categories for $N=2$ unitary minimal models}

Let $m_{i}$, $i=1, \dots, n$, be $n$ positive integers and
$V$ a vertex operator superalgebra in the class 
$\mathcal{C}_{m_{1};
\dots; m_{n}}$. Using  Corollary \ref{rat}, Proposition
\ref{2-7} and Theorem \ref{cep1} above, and Theorems 3.1 and 3.2
in \cite{H2}, which are  proved in 
 \cite{H2} using results
in \cite{HL1}--\cite{HL5} and \cite{H1}, we obtain the
following:

\begin{thm}
The conclusions of Theorems 4.1, 4.2, and 4.6
and the first conclusion of Corollary 4.7 in Section 4 of \cite{HM}
holds for $V$. \epf
\end{thm}

The notions of intertwining operator algebra 
in \cite{H2.5} (see also \cite{H4} and \cite{H6})
and $N=1$ superconformal intertwining 
operator superalgebra in \cite{HM} 
can be generalized easily to the following notion:

\begin{defn}
{\rm An {\it $N=2$ superconformal intertwining 
operator superalgebra} is an intertwining operator superalgebra
$W$
together with three elements $\tau^{+}$, $\tau^{-}$ and $\mu$
such that $(W^{e}, Y, \mathbf{1}, \tau^{+}, \tau^{-})$ is 
an $N=2$ superconformal vertex operator algebra.}
\end{defn}

Then we have:

\begin{thm}
The conclusion of Theorem 4.4 in \cite{HM} with ``$N=1$" repalced by ``$N=2$"
holds for $V$. \epf
\end{thm}

\noindent {\small \sc Department of Mathematics, Kerchof Hall, 
University of Virginia, Charlottesville, VA 22904-4137}

\noindent {\it and}

\noindent {\small \sc Department of Mathematics, Rutgers University,
110 Frelinghuysen Rd., Piscataway, NJ 08854-8019 (permanent address)}

\noindent {\em E-mail address}: yzhuang@math.rutgers.edu

\vskip 1em

\noindent {\small \sc Department of Mathematics, Rutgers University,
110 Frelinghuysen Rd., Piscataway, NJ 08854-8019}

\noindent {\em E-mail address}: amilas@math.rutgers.edu

\end{document}